\documentclass[12pt]{article}
\usepackage[dvipsnames]{xcolor}
\usepackage{amsmath}
\usepackage{setspace}
\usepackage{bm}
\usepackage{bbm}
\usepackage{amssymb}
\usepackage{indentfirst}
\usepackage[superscript]{cite}
\usepackage{geometry}
\usepackage{mathtools}
\usepackage{tocloft}

\DeclarePairedDelimiter\abs{\lvert}{\rvert}%
\DeclarePairedDelimiter{\floor}{\lfloor}{\rfloor}

\DeclareRobustCommand{\Stir}{\genfrac\{\}{0pt}{}}

\newcommand*{\citena}[1]{%
\begingroup
[\color{Green}
\romannumeral-`\x 
\setcitestyle{numbers}%
\cite{#1}%
\endgroup
]\ignorespacesafterend
}

\newcommand*{\citesup}[1]{%
\begingroup
\color{Green}
\cite{#1}%
\endgroup
\ignorespacesafterend
}

\newcommand*{\eqrefe}[1]{%
\begingroup
(\color{BrickRed}
\romannumeral-`\x 
\setcitestyle{numbers}%
\ref{eq:#1}%
\endgroup
)\ignorespacesafterend
}

\newcommand*{\secrefe}[1]{%
\begingroup
(\color{Aquamarine}
\romannumeral-`\x 
\setcitestyle{numbers}%
\ref{#1}%
\endgroup
)\ignorespacesafterend
}

\newtagform{Tags}[\textcolor{BrickRed}]{\color{Black}(}{)}

\usepackage[super]{natbib}
\setcitestyle{super}

\newcommand{\ii}{\bm{i}}

\begin{document}
\title{Lerch's $\Phi$ and the Polylogarithm at the Negative Integers}
\date{September 17, 2021}
\author{Jose Risomar Sousa}
\maketitle
\usetagform{Tags}

\begin{abstract}
At the negative integers, there is a simple relation between the Lerch $\Phi$ function and the polylogarithm. Starting from that relation and a formula for the polylogarithm at the negative integers known from the literature, we can deduce a simple closed formula for the Lerch $\Phi$ function at the negative integers, where the Stirling numbers of the second kind are not needed. Leveraging that finding, we also produce alternative formulae for the $k$-th derivatives of the cotangent and cosecant (ditto, tangent and secant), as simple functions of the negative polylogarithm and Lerch $\Phi$, respectively, which is evidence of the importance of these functions (they are less exotic than they seem). Lastly, we extend formulae for the Hurwitz zeta function only valid at the positive integers to the complex half-plane using this novelty.
\end{abstract}


\section{Introduction}
To start with some background, the Lerch transcendent $\Phi$ is a function defined by an infinite series,
\begin{equation} \nonumber
\Phi(e^{z},k,b)=\sum_{n=0}^{\infty}\frac{e^{z\,n}}{(n+b)^k} \text{,}
\end{equation}
\noindent that was introduced in 1887. It is a series that has the polylogarithm, the Hurwitz zeta and the Riemann zeta functions as particular cases.\\

Although not the main focus of this paper, it is important to note that if the real part of $z$ is negative, the $\Phi(e^{z},k,b)$ series always converges regardless of the sign of $k$ (assuming $k$ is real). Conversely, it never converges if $z$ has positive real part. Therefore, it only makes sense to talk about the analytic continuation of $\Phi(e^{z},k,b)$ in the $z$ plane, since the series does not converge when the real part of $z$ is positive. If the real part of $z$ is zero, however, we have a totally different scenario and analytic continuation is in the $k$ plane. The reason why only the real part of $z=x+\ii\,y$ matters for the convergence of such a series is the identity $e^{(x+\ii\,y)\,n}=e^{x\,n}\cdot e^{\ii\,y\,n}$ ($x$ and $y$ real). It implies that only the part $e^{x\,n}$ can explode out to infinity, as $e^{\ii\,y\,n}$ is a complex number where both the real and imaginary parts are less than one in modulus (in fact, $\abs{e^{\ii\,y\,n}}=1$). In the context of this paper, it is also important to make a distinction between the process of extending a formula from a narrower to a broader domain and analytic continuation. Section \secrefe{Lerch_AC} provides an example of the latter, while section \secrefe{Hurwitz_ext} illustrates the former.\\

As seen in reference \citena{Hurwitz}, it is possible to derive a formula for the Hurwitz zeta function at the positive integers with results from two previous papers that introduced new formulae for the generalized harmonic numbers and progressions, \citena{GHN} and \citena{CHP}, respectively.\\

To greatly summarize the reasoning presented in \citena{Hurwitz}, if $k$ is an integer greater than one then,
\begin{equation} \nonumber
\int_0^1 u^k(1-\cos{2\pi\,n\,u})\cot{\pi u}\,du \sim -\frac{H(n)}{\pi}-\frac{k!}{\pi}\sum_{j=1}^{\floor{(k-1)/2}}\frac{(-1)^j(2\pi)^{-2j}\zeta(2j+1)}{(k-2j)!} \text{,}
\end{equation}
\noindent which implies the following approximation,
\begin{equation} \nonumber
\int_0^1 (u^k-u)(1-\cos{2\pi\,n\,u})\cot{\pi u}\,du \sim -\frac{k!}{\pi}\sum_{j=1}^{\floor{(k-1)/2}}\frac{(-1)^j(2\pi)^{-2j}\zeta(2j+1)}{(k-2j)!}=\int_0^1 (u^k-u)\cot{\pi u}\,du \text{,}
\end{equation}
\noindent with the equality on the right only valid for positive integer $k$. And this approximation in turn justifies the formula\footnote{The formulae that can be derived with this method are not unique and the one shown may be the simplest.} shown next.\\

For every integer $k$ greater than one and every non-integer complex $b$,
\begin{multline} \label{eq:Hurwitz_form1}
\zeta(k,b)=\frac{1}{2b^k}+\frac{(2\pi\ii)^k}{4}\left(\frac{\mathrm{Li}_{-k+1}\left(e^{-2\pi\ii\,b}\right)}{(k-1)!}+e^{-2\pi\ii\,b}\sum_{j=1}^{k}\frac{\delta_{1j}+\mathrm{Li}_{-j+1}(e^{-2\pi\ii\,b})}{(j-1)!(k-j)!}\right)\\-\frac{\ii(2\pi\ii)^k}{2}\int_{0}^{1}\sum_{j=1}^{k}\frac{\left(\delta_{1j}+\mathrm{Li}_{-j+1}(e^{-2\pi\ii\,b})\right)\left(u^{k-j}e^{-2\pi\ii\,b\,u}-e^{-2\pi\ii\,b}\right)}{(j-1)!(k-j)!}\cot{\pi u}\,du 
\end{multline}
\indent The discovery of a new, possibly first, closed formula for the Lerch transcendent function at the negative integers was made possible through the analysis of the above formula.\\

The big breakthrough is really in the following straightforward identity, which only works for the analytic continuation of the Lerch $\Phi$ and the polylogarithm functions at the negative integers, and makes it possible to obtain the former as a sum of the latter,
\begin{equation} \label{eq:Soma_Polylog}
\Phi(e^{z},-k,b+1)=e^{-z}\,\sum_{j=0}^{k}\binom{k}{j}\mathrm{Li}_{-j}\left(e^{z}\right)b^{k-j}
\end{equation}
\indent In fact, this identity also applies to the sum of Lerches,
\begin{equation} \label{eq:Soma_Lerch}
\Phi(e^{z},-k,u+v)=\sum_{j=0}^{k}\binom{k}{j}\Phi\left(e^{z},-j,v\right)u^{k-j}
\end{equation}
\indent Since formulae \eqrefe{Soma_Polylog} and \eqrefe{Soma_Lerch} break down when $e^z=1$, from two known facts from the literature, namely, a recurrence for the Bernoulli polynomials and a relation between these and the Hurwitz zeta function,
\begin{equation} \nonumber
B_k(u+v)=\sum_{j=0}^{k}\binom{k}{j}B_j(v)\,u^{k-j} \text{ and } B_j(v)=-j\,\zeta(1-j,v) \text{,}
\end{equation}
\noindent one can conclude that,
\begin{equation} \label{eq:Soma_Hurwz}
\zeta(-k,u+v)=-\frac{u^{k+1}}{k+1}+\sum_{j=0}^{k}\binom{k}{j}\zeta(-j,v)\,u^{k-j} \text{,} 
\end{equation}
\noindent and from this recurrence, a natural expression relating the Hurwitz and the Riemann zeta functions can be obtained, which completes the picture,
\begin{equation} \label{eq:Soma_Hurwz_part}
\zeta(-k,b+1)=-\frac{b^{k+1}}{k+1}+\sum_{j=0}^{k}\binom{k}{j}\zeta(-j)\,b^{k-j}
\end{equation}
\indent Relation \eqrefe{Soma_Polylog} may be the counterpart to the functional equation that relates the Lerch $\Phi$ to the Hurwitz zeta at the positive integers, as demonstrated in a previous paper\citesup{Lerch}. This new relation therefore presupposes the existence of another relation between the polylogarithm and the zeta function at the negative integers.

\section{Stirling numbers of the second kind}
The existing formula for the polylogarithm, $\mathrm{Li}_{-j+1}(z)$, available in the literature, makes use of the Stirling numbers of the second kind\footnote{The $S(j,q)$ in the curly brackets.}. If $j$ is a positive integer then,
\begin{equation} \label{eq:Polylog_lit}
\mathrm{Li}_{-j+1}(z)=\sum_{q=1}^{j}(q-1)!\Stir{j}{q}\left(\frac{z}{1-z}\right)^q
\end{equation}
\indent Let us see how the Lerch $\Phi$ can be obtained from it.  First, formula \eqrefe{Polylog_lit} is transformed to resemble \eqrefe{Hurwitz_form1},
\begin{equation} \nonumber
\mathrm{Li}_{-j+1}(e^{-2\pi\ii\,z})=\sum_{q=1}^{j}(q-1)!\Stir{j}{q}\left(-\frac{1+\ii\cot{\pi\,z}}{2}\right)^q
\end{equation}
\indent Now, going back to equation \eqrefe{Soma_Polylog},
\begin{equation} \nonumber
\Phi(e^{-2\pi\ii\,z},-k+1,b+1)=e^{2\pi\ii\,z}\,\sum_{j=1}^{k}\binom{k-1}{j-1}\mathrm{Li}_{-j+1}\left(e^{-2\pi\ii\,z}\right)b^{k-j} \Rightarrow
\end{equation}
\begin{equation} \nonumber
\Phi(e^{-2\pi\ii\,z},-k+1,b+1)=e^{2\pi\ii\,z}\,\sum_{j=1}^{k}\binom{k-1}{j-1}\,b^{k-j}\sum_{q=1}^{j}(q-1)!\Stir{j}{q}\left(-\frac{1+\ii \cot{\pi\,z}}{2}\right)^q \Rightarrow
\end{equation}
\begin{equation} \label{eq:Lerch_form}
\Phi(e^{-2\pi\ii\,z},-k+1,b+1)=e^{2\pi\ii\,z}\,\sum_{q=1}^{k}(q-1)!\left(-\frac{1+\ii \cot{\pi\,z}}{2}\right)^q\sum_{j=q}^{k}\binom{k-1}{j-1}\Stir{j}{q}\,b^{k-j}
\end{equation}\\
\indent At this point we need to inspect the inner sum from equation \eqrefe{Lerch_form} and see if it is possible to rewrite it. An expression for a similar sum exists in the literature, \begin{equation} \label{eq:Stir_lit}
\sum_{j=q}^{k}\binom{k}{j}\Stir{j}{q}=\Stir{k+1}{q+1} \text{,}
\end{equation}
\noindent but this new one is much more complicated.

\section{Binomial formula for Stirling numbers}
Let us try and rewrite the below finite sum,
\begin{equation} \label{eq:Stir_soma}
\sum_{j=q}^{k}\binom{k-1}{j-1}\Stir{j}{q}\,b^{k-j}
\end{equation}
\indent First off, the literature provides us with the below relation,
\begin{equation} \nonumber
\sum_{j=q}^{\infty}\Stir{j}{q}\frac{x^{j}}{j!}=\frac{(e^x-1)^q}{q!} \text{,}
\end{equation}
\noindent therefore,
\begin{equation} \nonumber
\sum_{k=q}^{\infty}\sum_{j=q}^{k}\Stir{j}{q}\frac{x^{j}}{j!}\frac{y^{k-j}}{(k-j)!}=\frac{e^y(e^x-1)^q}{q!} \text{,}
\end{equation}\\
\noindent and differentiating with respect to $x$,
\begin{equation} \nonumber
\sum_{k=q}^{\infty}\sum_{j=q}^{k}\Stir{j}{q}\frac{x^{j-1}}{(j-1)!}\frac{y^{k-j}}{(k-j)!}=\frac{e^y e^x(e^x-1)^{q-1}}{(q-1)!} \text{}
\end{equation}
\indent Now, making $x=b\,z$ and $y=z$, one has,
\begin{equation} \nonumber
\sum_{k=q}^{\infty}z^{k-1}\sum_{j=q}^{k}\Stir{j}{q}\frac{b^{j-1}}{(j-1)!(k-j)!}=\frac{e^{(b+1)z}(e^{b\,z}-1)^{q-1}}{(q-1)!}=\sum_{j=0}^{q-1}\frac{(-1)^{q-1-j}\,e^{(b+1+j\,b)z}}{j!(q-1-j)!} \text{,}
\end{equation}
\noindent where the rightmost expression stems from the Newton binomial. Hence, differentiating $k-1$ times with respect to $z$, one concludes that,
\begin{equation} \label{eq:Stir_final}
\sum_{j=q}^{k}\binom{k-1}{j-1}\Stir{j}{q}\,b^{k-j}=\sum_{j=1}^{q}\frac{(-1)^{q-j}(j+b)^{k-1}}{(j-1)!(q-j)!} \text{}
\end{equation}\\
\indent To obtain a more appropriate version of this formula, it can be integrated as below,
\begin{equation} \nonumber
\int_{0}^{b}\sum_{j=q}^{k}\frac{x^{j-1}}{(j-1)!(k-j)!}\Stir{j}{q}\,dx=\frac{1}{(k-1)!}\int_{0}^{b}\sum_{j=1}^{q}\frac{(-1)^{q-j}(j\,x+1)^{k-1}}{(j-1)!(q-j)!}\,dx \text{,}
\end{equation}
\noindent which gives the neater expression,
\begin{equation} \label{eq:Stir_final2}
\sum_{j=q}^{k}\binom{k}{j}\Stir{j}{q}\,b^{k-j}=\sum_{j=0}^{q}\frac{(-1)^{q-j}(j+b)^{k}}{j!(q-j)!} \text{,}
\end{equation}
\noindent which holds for every non-negative integer $q$ and every $b$.\\
 
\indent Though it is not going to be used here, another pattern similar to the binomial theorem emerges in the factorial power of the sum of two numbers, $x$ and $y$,
\begin{equation} \nonumber
(x+y)^{(k)}=\sum_{j=0}^{k}\binom{k}{j}x^{(j)}\,y^{(k-j)} \text{, where } x^{(j)}=\frac{x!}{(x-j)!}
\end{equation}

\section{Lerch's $\Phi$ at the negative integers} \label{Lerch_AC}
When equations \eqrefe{Soma_Polylog}, \eqrefe{Lerch_form} and \eqrefe{Stir_final} are combined, the result is the following formula,
\begin{equation} \label{eq:Lerch_complex_final}
\Phi(e^{-2\pi\ii\,z},-k+1,b+1)=e^{2\pi\ii\,z}\,\sum_{q=1}^{k}\left(\frac{1+\ii \cot{\pi\,z}}{2}\right)^q\sum_{j=1}^{q}\binom{q-1}{j-1}(-1)^j\,(j+b)^{k-1} \text{}
\end{equation}
\indent For integer $z$, the formula is not defined as the cotangent is infinity, so we can not extract the Hurwitz zeta at the negative integers from it. But from the relation,
\begin{equation} \nonumber
\Phi(e^{-2\pi\ii\,z},-k+1,1)=e^{2\pi\ii\,z}\,\mathrm{Li}_{-k+1}\left(e^{-2\pi\ii\,z}\right) \text{,}
\end{equation}
\noindent the polylogarithm can be derived, which however is just a rewrite of equation \eqrefe{Polylog_lit} with an expression for $S(k,q)$ known from the literature, which nonetheless confirms equation \eqrefe{Lerch_complex_final},
\begin{equation} \label{eq:Polylog_complex_final}
\mathrm{Li}_{-k+1}\left(e^{-2\pi\ii\,z}\right)=\sum_{q=1}^{k}\left(\frac{1+\ii \cot{\pi\,z}}{2}\right)^q\sum_{j=1}^{q}\binom{q-1}{j-1}(-1)^j\,j^{k-1}
\end{equation}
\indent Looking at formulae \eqrefe{Lerch_complex_final} and \eqrefe{Polylog_complex_final} now, it might look simple to go from the latter straight to the former without having to solve \eqrefe{Stir_soma}, but that is misleading.\\

Finally, \eqrefe{Lerch_complex_final} can be turned into a simpler form, which holds for every non-negative integer $k$,
\begin{equation} \label{eq:Lerch_final}
\Phi(z,-k,\,b)=-\frac{1}{z-1}\sum_{q=0}^{k}\left(\frac{z}{z-1}\right)^q\sum_{j=0}^{q}\binom{q}{j}(-1)^j\,(j+b)^{k}
\end{equation}\\
\indent The above expression gives the analytic continuation of the Lerch $\Phi$ at the non-positive integers $-k$ (since it holds for all $z$, except $z=1$). It is interesting to note how much simpler the formula of the Lerch $\Phi$ at the negative integers is than the formula at the positive integers from \citena{Lerch}. And also how strikingly similar it is to the power series for the Lerch $\Phi$ available in the literature, which holds for all $k$ and $z$ with $\Re(z)<1/2$,
\begin{equation} \nonumber
\Phi(z,k,\,b)=-\frac{1}{z-1}\sum_{q=0}^{\infty}\left(\frac{z}{z-1}\right)^q\sum_{j=0}^{q}\binom{q}{j}(-1)^j\,(j+b)^{-k}
\end{equation}

\section{Derivatives of trigonometric functions}
In his paper \textit{On the Hurwitz function for rational arguments}\citesup{Adam}, Victor Adamchik provides the first ever formula for the intricate patterns of the $k$-th derivatives of the cotangent. It looks like this,
\begin{equation} \nonumber
\frac{\text{d}^k(\cot{a\,x})}{\text{d}\,x^k}=(2\,\ii\,a)^k(-\ii+\cot{a\,x})\sum_{q=1}^{k}q!\Stir{k}{q}\left(-\frac{1-\ii \cot{a\,x}}{2}\right)^q
\end{equation}
\indent It is possible to express this formula as a simple function of the polylogarithm. First, we rewrite it as,
\begin{equation} \label{eq:Cot}
\frac{\text{d}^k(\cot{a\,x})}{\text{d}\,x^k}=(2\,\ii\,a)^k(-\ii+\cot{a\,x})\sum_{q=1}^{k}\left(\frac{1-\ii \cot{a\,x}}{2}\right)^q\sum_{j=1}^{q}\frac{q!\,(-1)^j\, j^{k-1}}{(j-1)!(q-j)!} \text{,}
\end{equation}
\noindent where $S(k,q)$ was replaced by an equivalent formula,
\begin{equation} \label{eq:Stir_form_vari}
\Stir{k}{q}=(-1)^q\,\sum_{j=1}^{q}\frac{(-1)^j\, j^{k-1}}{(j-1)!(q-j)!} \text{,}
\end{equation}
\noindent that stems from equations \eqrefe{Stir_lit} and \eqrefe{Stir_final2}.\\

Secondly, we note how similar it looks to the polylogarithm from \eqrefe{Polylog_complex_final},
\begin{equation} \nonumber
\mathrm{Li}_{-k+1}\left(e^{2\,\ii\,a\,x}\right)=\sum_{q=1}^{k}\left(\frac{1-\ii \cot{a\,x}}{2}\right)^q\sum_{j=1}^{q}\frac{(q-1)!\,(-1)^j\, j^{k-1}}{(j-1)!(q-j)!}
\end{equation}
\indent If the above polylog is differentiated once with respect to $x$ and transformed, an alternative expression is obtained for the polylogarithm of order $k$,
\begin{equation} \label{eq:Poly_Transf}
\mathrm{Li}_{-k}\left(e^{2\,\ii\,a\,x}\right)=\frac{1}{1-e^{2\,\ii\,a\,x}}\sum_{q=1}^{k}\left(\frac{1-\ii \cot{a\,x}}{2}\right)^q\sum_{j=1}^{q}\frac{q!\,(-1)^j\, j^{k-1}}{(j-1)!(q-j)!} \text{,}
\end{equation}
\noindent which, however, is not exactly equal to form \eqrefe{Polylog_complex_final}. That stems from a property of polylogs, that when differentiated they yield the next order polylog.\\

Finally, comparing the two expressions, \eqrefe{Cot} and \eqrefe{Poly_Transf}, we conclude that,
\begin{equation} \label{eq:Cot_final}
\frac{\text{d}^k(\cot{a\,x})}{\text{d}\,x^k}=-\ii\,\delta_{0\,k}-2\,\ii(2\,\ii\,a)^k\,\mathrm{Li}_{-k}\left(e^{2\,\ii\,a\,x}\right) \text{, where } \delta_{0\,k}=1 \text{ iff }k=0
\end{equation}
\indent To obtain the cosecant, we can resort to a simple logic,
\begin{equation} \nonumber
\frac{\cos{a\,x}+\ii\sin{a\,x}}{\sin{a\,x}}=\frac{e^{\ii\,a\,x}}{\sin{a\,x}}=\ii+\cot{a\,x} \Rightarrow \frac{1}{\sin{a\,x}}=e^{-\ii\,a\,x}(\ii+\cot{a\,x}) \text{,}
\end{equation}\\
\noindent and then the Leibniz rule for the derivative of a product of two functions leads to,
\begin{equation} \nonumber
\frac{\text{d}^k}{\text{d}\,x^k}\left(\frac{1}{\sin{a\,x}}\right)=-2\,\ii\,e^{-\ii\,a\,x}\sum_{q=0}^{k}\binom{k}{q}(2\,\ii\,a)^q\,\mathrm{Li}_{-q}\left(e^{2\,\ii\,a\,x}\right)(-\ii\,a)^{k-q}
\end{equation}
\indent Lastly, formula \eqrefe{Soma_Polylog} allows the above expression to be rewritten as,
\begin{equation} \label{eq:Cosec}
\frac{\text{d}^k}{\text{d}\,x^k}\left(\frac{1}{\sin{a\,x}}\right)=-2\,\ii(2\,\ii\,a)^k\,e^{\ii\,a\,x}\,\Phi\left(e^{2\,\ii\,a\,x},-k,\frac{1}{2}\right) \text{,}
\end{equation}
\noindent which holds for every non-negative integer $k$.

\subsection{Tangent and secant}
To be able to obtain the tangent and secant, first we need to produce a formula for the cotangent and cosecant of a translated arc. Adamchik's formula\citesup{Adam} becomes, 
\begin{equation} \nonumber
\frac{\text{d}^k\left(\cot{(a\,x+b)}\right)}{\text{d}\,x^k}=(2\,\ii\,a)^k\left(-\ii+\cot{(a\,x+b)}\right)\sum_{q=1}^{k}\left(\frac{1-\ii \cot{(a\,x+b)}}{2}\right)^q\sum_{j=1}^{q}\frac{q!\,(-1)^j\,j^{k-1}}{(j-1)!(q-j)!}
\end{equation}
\indent The polylog formula then changes to,
\begin{equation} \nonumber
\mathrm{Li}_{-k}\left(e^{2\,\ii(a\,x+b)}\right)=\frac{1}{1-e^{2\,\ii(a\,x+b)}}\sum_{q=1}^{k}\left(\frac{1-\ii \cot{(a\,x+b)}}{2}\right)^q\sum_{j=1}^{q}\frac{q!\,(-1)^j\,j^{k-1}}{(j-1)!(q-j)!} \text{,}
\end{equation}
\noindent and then the final formula is not too different from the simple case,
\begin{equation} \label{Cot_trans}
\frac{\text{d}^k\left(\cot{(a\,x+b)}\right)}{\text{d}\,x^k}=-\ii\,\delta_{0\,k}-2\,\ii(2\,\ii\,a)^k\,\mathrm{Li}_{-k}\left(e^{2\,\ii(a\,x+b)}\right) \text{, where } \delta_{0\,k}=1 \text{ iff }k=0
\end{equation}\\
\indent Similarly, the cosecant of a translated arc is,
\begin{equation} \label{Cosec_trans}
\frac{\text{d}^k}{\text{d}\,x^k}\left(\frac{1}{\sin{(a\,x+b)}}\right)=-2\,\ii(2\,\ii\,a)^k\,e^{\ii(a\,x+b)}\,\Phi\left(e^{2\,\ii(a\,x+b)},-k,\frac{1}{2}\right) \text{}
\end{equation}\\
\indent Finally, to obtain the tangent and secant, we just need to set $b$ to $\pi/2$. And since the formulae for the translated arc are not very different from when $b=0$, for the tangent one has,
\begin{equation} \label{Tan_trans}
\frac{\text{d}^k\left(\tan{(a\,x+b)}\right)}{\text{d}\,x^k}=\ii\,\delta_{0\,k}+2\,\ii(2\,\ii\,a)^k\,\mathrm{Li}_{-k}\left(-e^{2\,\ii(a\,x+b)}\right) \text{, where } \delta_{0\,k}=1 \text{ iff }k=0 \text{,}
\end{equation}
\noindent and for the secant,
\begin{equation} \label{eq:Sec_trans}
\frac{\text{d}^k}{\text{d}\,x^k}\left(\frac{1}{\cos{(a\,x+b)}}\right)=2\,(2\,\ii\,a)^k\,e^{\ii(a\,x+b)}\,\Phi\left(-e^{2\,\ii(a\,x+b)},-k,\frac{1}{2}\right) \text{}
\end{equation}\\
\indent It is surprising that these derivatives can be expressed by means of negative Lerch and polylogs. For example, the negative polylog is known to yield the derivatives of a simple exponential function at a point, but not the derivative itself,
\begin{equation} \nonumber
\frac{\text{d}^k}{\text{d}\,x^k}\left(\frac{x}{e^{a\,x+b}-1}\right)\bigg|_{x=0}=-k\left(\delta_{1\,k}+\mathrm{Li}_{-k+1}\left(e^{b}\right)\right)a^{k-1}
\end{equation}

\section{The extended Hurwitz zeta formula} \label{Hurwitz_ext}
The domain of the Hurwitz zeta formula from \eqrefe{Hurwitz_form1} can be extended from the integers greater than one to the complex numbers $k$ with real part greater than one. From the relation \eqrefe{Soma_Polylog} one obtains, 
\begin{equation} \nonumber
\sum_{j=1}^{k}\frac{\mathrm{Li}_{-j+1}\left(e^{-2\pi\ii\,b}\right)u^{k-j}}{(j-1)!(k-j)!}=\frac{e^{-2\pi\ii\,b}}{(k-1)!}\Phi(e^{-2\pi\ii\,b},-k+1,u+1)\text{,}
\end{equation} 
\noindent which replaced into \eqrefe{Hurwitz_form1} gives the below,
\begin{multline} \nonumber
\zeta(k,b)=\frac{1}{2b^k}+\frac{(2\pi\ii)^k}{4\,\Gamma(k)}\left(e^{-2\pi\ii\,b}+e^{-4\pi\ii\,b}\,\Phi(e^{-2\pi\ii\,b},-k+1,2)+\mathrm{Li}_{-k+1}\left(e^{-2\pi\ii\,b}\right)\right)\\-\frac{\ii(2\pi\ii)^k}{2\,\Gamma(k)}\int_{0}^{1}\Big(u^{k-1}e^{-2\pi\ii\,b\,u}-e^{-2\pi\ii\,b}\\ +e^{-2\pi\ii\,b(u+1)}\,\Phi(e^{-2\pi\ii\,b},-k+1,u+1)-e^{-4\pi\ii\,b}\,\Phi(e^{-2\pi\ii\,b},-k+1,2)\Big)\cot{\pi u}\,du \text{}
\end{multline}\\
\indent The above formula can be simplified further with the identity,
\begin{equation} \label{eq:Lerch_neg_id}
u^{k-1}+e^{-2\pi\ii\,b}\,\Phi(e^{-2\pi\ii\,b},-k+1,u+1)=\Phi(e^{-2\pi\ii\,b},-k+1,u) \text{}
\end{equation}
\indent That leads to the simpler form below, which holds when $\Re{(k)}>1$ and $b$ is not integer,
\begin{multline} \label{eq:Hurwitz}
\zeta(k,b)=\frac{1}{2b^k}+\frac{(2\pi\ii)^k}{2\,\Gamma(k)}\,\mathrm{Li}_{-k+1}\left(e^{-2\pi\ii\,b}\right)\\-\frac{\ii(2\pi\ii)^k}{2\,\Gamma(k)}\int_{0}^{1}\Big(e^{-2\pi\ii\,b\,u}\,\Phi(e^{-2\pi\ii\,b},-k+1,u)-\mathrm{Li}_{-k+1}\left(e^{-2\pi\ii\,b}\right)\Big)\cot{\pi u}\,du \text{}
\end{multline}\\
\indent Since the $\zeta(k,b)$ series converges for $\Re{(k)}>1$, this is not an analytic continuation, it is just a way to extend the formula beyond the integers greater than one.

\subsection{Hurwitz zeta formula rewritten} 
Now that fomula \eqrefe{Lerch_complex_final} is known, the Hurwitz zeta formula at the integers greater than one from \eqrefe{Hurwitz_form1} can be rewritten with only references to elementary functions.\\

The pattern of the formula now becomes more apparent than in \eqrefe{Hurwitz_form1}, as two terms of the integrand previously not included into the summation symbol can now be moved under it, which is accomplished by means of the binomial coefficient,
\begingroup
\small
\begin{multline} \nonumber
\zeta(k,b)=\frac{1}{2b^k}+\frac{(2\pi\ii)^k}{4(k-1)!}\sum_{q=0}^{k}\left(\frac{1+\ii \cot{\pi\,b}}{2}\right)^q\sum_{j=0}^{q}\binom{q-1}{j-1}(-1)^j\left(e^{-2\pi\ii\,b}(j+1)^{k-1}+j^{k-1}\right)\\ -\frac{\ii(2\pi\ii)^k}{2(k-1)!}\int_{0}^{1}\sum_{q=0}^{k}\left(\frac{1+\ii \cot{\pi\,b}}{2}\right)^q\sum_{j=0}^{q}\binom{q-1}{j-1}(-1)^j\left(e^{-2\pi\ii\,b\,u}(j+u)^{k-1}-e^{-2\pi\ii\,b}(j+1)^{k-1}\right)\cot{\pi u}\,du 
\end{multline} 
\endgroup\\
\indent One of the advantages of this new formula is the fact it allows one to get rid of its non-real parts more easily, though the resulting formula is inevitably more complicated.

\subsection{When the parameter $b$ is a half-integer}
The below result stems from $\cot{\pi\,b}=0$ and $e^{-2\pi\ii\,b}=-1$ when $b$ is a half-integer,
\begin{multline} \label{eq:Hurwitz_half}
\zeta(k,b)=\frac{1}{2b^k}-\frac{(2\pi\ii)^k}{4(k-1)!}\sum_{q=0}^{k}\left(\frac{1}{2}\right)^q\sum_{j=0}^{q}\binom{q-1}{j-1}(-1)^j\left((j+1)^{k-1}-j^{k-1}\right)\\ -\frac{\ii(2\pi\ii)^k}{2(k-1)!}\int_{0}^{1}\sum_{q=0}^{k}\left(\frac{1}{2}\right)^q\sum_{j=0}^{q}\binom{q-1}{j-1}(-1)^j\left(e^{-2\pi\ii\,b\,u}(j+u)^{k-1}+(j+1)^{k-1}\right)\cot{\pi u}\,du 
\end{multline} 

\section{A new formula for the Hurwitz zeta}
In \citena{Hurwitz}, we had created a generating function for the Hurwitz zeta function, $f(x)$. When $b$ is not a half-integer or integer, the expression is,
\begin{multline} \label{eq:f(x)}
f(x)=\sum _{k=2}^{\infty}x^k\zeta(k,b)=-\frac{x^2}{2b(x-b)}-\frac{1}{2\sin{\pi b}}\frac{\pi x\sin{\pi x}}{\sin{\pi(x-b)}}\\-\pi x\int _0^1\left(\frac{\sin{2\pi u(x-b)}}{\sin{2\pi(x-b)}}-\frac{\sin{2\pi b u}}{\sin{2\pi b}}\right)\cot{\pi u}\,du
\end{multline}
\indent The $k$-th derivative of $f(x)$ yields the Hurwitz zeta function of order $k$,
\begin{equation} \nonumber
\zeta(k,b)=\frac{f^{(k)}(0)}{k!} \text{}
\end{equation}
\indent And now that we know how to differentiate the cosecant successively, it is possible to produce an explicit formula from $f(x)$, again through the Leibniz rule. However, to make this process simpler, we resort to two artifices. First, to get rid of the extra $x$ factor in the integral, we divide $f(x)$ by $x$ and take the $(k-1)$-th derivative instead of the $k$-th. Second, to avoid the complications of differentiating the sine, we replace it with an equivalent sum of exponential functions.\\

The first and second parts of \eqrefe{f(x)} are straightforward, they coincide with the terms outside of the integral from \eqrefe{Hurwitz}, that is,
\begingroup
\small
\begin{equation} \nonumber
\frac{1}{k!}\frac{\text{d}^k}{\text{d}\,x^k}\left(-\frac{1}{2\sin{\pi b}}\frac{\pi x\sin{\pi x}}{\sin{\pi(x-b)}}\right)\bigg|_{x=0}=\frac{(2\pi\ii)^k}{2\,(k-1)!}\,\mathrm{Li}_{-k+1}\left(e^{-2\pi\ii\,b}\right)
\end{equation}
\endgroup\\
\indent The same is not true for the third part of \eqrefe{f(x)}, since $f(x)$ was created using a different process than \eqrefe{Hurwitz} (see \citena{Hurwitz} for details). The integrals evaluate to the same number, but the integrands are not the same.\\

After all is put together, the final formula holds for any integer $k$ greater than one and any $b$ that is not an integer or a half-integer,
\begingroup
\small
\begin{multline} \label{eq:Hurwitz_2}
\zeta(k,b)=\frac{1}{2b^k}+\frac{(2\pi\ii)^k}{2\,(k-1)!}\,\mathrm{Li}_{-k+1}\left(e^{-2\pi\ii\,b}\right)+ \\ -\frac{\ii(2\pi\ii)^k\,e^{-2\pi\ii\,b}}{4}\int_{0}^{1}\sum_{j=1}^{k}\frac{2^j u^{k-j}\left(e^{-2\pi\ii\,b\,u}-(-1)^{k-j}e^{2\pi\ii\,b\,u}\right)}{(j-1)!(k-j)!}\Phi\left(e^{-4\pi\ii\,b},-j+1,\frac{1}{2}\right)\cot{\pi u}\,du 
\end{multline} 
\endgroup\\
\indent Finally, using the relation \eqrefe{Soma_Lerch}, the formula can be extended to $\Re(k)>1$,
\begingroup
\small
\begin{multline} \nonumber
\zeta(k,b)=\frac{1}{2b^k}+\frac{(2\pi\ii)^k}{2\,\Gamma(k)}\,\mathrm{Li}_{-k+1}\left(e^{-2\pi\ii\,b}\right)+ \\ -\frac{\ii(4\pi\ii)^k\,e^{-2\pi\ii\,b}}{4\,\Gamma(k)}\int_{0}^{1}\left(e^{-2\pi\ii\,b\,u}\,\Phi\left(e^{-4\pi\ii\,b},-k+1,\frac{u+1}{2}\right)-e^{2\pi\ii\,b\,u}\,\Phi\left(e^{-4\pi\ii\,b},-k+1,\frac{-u+1}{2}\right)\right)\cot{\pi u}\,du 
\end{multline}\\
\endgroup


\end{document}